# EXTRA HEADS AND INVARIANT ALLOCATIONS

By Alexander E. Holroyd[1] and Yuval Peres[2]

*University of British Columbia and University of California, Berkeley*

Let $\Pi$ be an ergodic simple point process on $\mathbb{R}^d$ and let $\Pi^*$ be its Palm version. Thorisson [*Ann. Probab.* **24** (1996) 2057–2064] proved that there exists a *shift coupling* of $\Pi$ and $\Pi^*$; that is, one can select a (random) point $Y$ of $\Pi$ such that translating $\Pi$ by $-Y$ yields a configuration whose law is that of $\Pi^*$. We construct shift couplings in which $Y$ and $\Pi^*$ are functions of $\Pi$, and prove that there is no shift coupling in which $\Pi$ is a function of $\Pi^*$. The key ingredient is a deterministic translation-invariant rule to allocate sets of equal volume (forming a partition of $\mathbb{R}^d$) to the points of $\Pi$. The construction is based on the Gale–Shapley stable marriage algorithm [*Amer. Math. Monthly* **69** (1962) 9–15]. Next, let $\Gamma$ be an ergodic random element of $\{0,1\}^{\mathbb{Z}^d}$ and let $\Gamma^*$ be $\Gamma$ conditioned on $\Gamma(0) = 1$. A shift coupling $X$ of $\Gamma$ and $\Gamma^*$ is called an *extra head scheme*. We show that there exists an extra head scheme which is a function of $\Gamma$ if and only if the marginal $\mathbf{E}[\Gamma(0)]$ is the reciprocal of an integer. When the law of $\Gamma$ is product measure and $d \geq 3$, we prove that there exists an extra head scheme $X$ satisfying $\mathbf{E} \exp c \|X\|^d < \infty$; this answers a question of Holroyd and Liggett [*Ann. Probab.* **29** (2001) 1405–1425].

**1. Introduction.** Let $\Pi$ be a translation-invariant ergodic simple point process of unit intensity on $\mathbb{R}^d$, with law $\Lambda$. Let $\Pi^*$ be the Palm version of $\Pi$, with law $\Lambda^*$. (Recall that if $\Pi$ is a Poisson process, $\Pi^*$ is a Poisson process with an added point at the origin.) We call elements of $\mathbb{R}^d$ *sites* and we call integer-valued Borel measures on $\mathbb{R}^d$ *configurations* (so $\Pi$ and $\Pi^*$ are random configurations). For a configuration $\pi$ and a site $y$ we write $T^{-y}\pi$ for the translated configuration given by $(T^{-y})\pi(\cdot) = \pi(\cdot + y)$. A (*continuum*)

Received June 2003; revised September 2003.

[1]Supported in part by NSF Grants DMS-00-72398 and CCR-01-21555, by an NSERC Discovery Grant and by the Center for Pure and Applied Mathematics at UC Berkeley.

[2]Supported in part by NSF Grants DMS-01-04073, DMS-02-44479, CCR-01-21555 and by a Miller Professorship at UC Berkeley.

*extra head scheme* for $\Pi$ is an $\mathbb{R}^d$-valued random variable $Y$ such that the point process $T^{-Y}\Pi$ has law $\Lambda^*$. Thorisson [13] proved (in a more general setting) that for any $\Pi$ as above, there exists a continuum extra head scheme. We may regard an extra head scheme as a *shift-coupling*, that is, a coupling $(\Pi, \Pi^*, Y)$ in which $\Pi, \Pi^*$ have respective laws $\Lambda, \Lambda^*$, and $\Pi^* = T^{-Y}\Pi$ almost surely. A *nonrandomized* extra head scheme is a shift coupling in which $Y$ (and therefore $\Pi^*$) is almost surely a function of $\Pi$. We shall prove the following.

THEOREM 1. *For any $d \geq 1$ and any translation-invariant ergodic simple point process $\Pi$ in $\mathbb{R}^d$, there exists a nonrandomized extra head scheme.*

Liggett [8] proved Theorem 1 in the case $d = 1$. In contrast, we have the following.

PROPOSITION 2. *Let $d \geq 1$ and let $\Pi$ be any ergodic translation-invariant simple point process on $\mathbb{R}^d$. For any shift coupling of $\Pi, \Pi^*$ where $\Pi = T^Y \Pi^*$, the translation variable $Y$ cannot be a function of $\Pi^*$.*

Given that extra head schemes exist, it is natural to ask how to construct an extra head scheme $Y$ from the configuration $\Pi$. The existence proof in [13] gives little clue how to do this; on the other hand, in [8], an explicit construction for a nonrandomized extra head scheme is given for $d = 1$. Our proof of Theorem 1 will be based on the following construction. The *support* of $\Pi$ is the random set $[\Pi] = \{x \in \mathbb{R}^d : \Pi(\{x\}) = 1\}$. A *balanced allocation rule* for $\Pi$ is a measurable function $\Psi_\Pi : \mathbb{R}^d \to [\Pi]$, defined from $\Pi$ in a deterministic, translation-invariant way, such that $\Psi_\Pi^{-1}(y)$ has Lebesgue measure 1 for each $y \in [\Pi]$. (We shall give a more careful definition later.) From a balanced allocation rule $\Psi$, we shall obtain a nonrandomized extra head scheme by taking $Y = \Psi_\Pi(0)$. We shall construct a balanced allocation rule using an approach based on the Gale–Shapley stable marriage algorithm [2]. The resulting $\Psi_\Pi$ is illustrated in Figure 1. Its properties are studied in detail in [4]. Related questions involving stable matchings of point processes were studied in [6].

Consider now the following discrete setting. Let $\mu$ be a translation-invariant ergodic measure on the product $\sigma$-algebra of $\{0,1\}^{\mathbb{Z}^d}$. We call elements of $\mathbb{Z}^d$ *sites* and elements of $\{0,1\}^{\mathbb{Z}^d}$ *configurations*. Let $\Gamma$ be a random configuration with law $\mu$. We say that a site $x$ is *occupied* if $\Gamma(x) = 1$ and *unoccupied* if $\Gamma(x) = 0$. Let $p$ be the marginal probability that the origin is occupied, and assume $p \in (0,1)$. Let $\mu^*$ be the conditional law of $\Gamma$ given that the origin is occupied. For a site $z$ and a configuration $\gamma$ we denote by $T^{-z}\gamma$ the translated configuration given by $(T^{-z}\gamma)(y) = \gamma(y + z)$. A *(discrete) extra head scheme* for $\Gamma$ is a $\mathbb{Z}^d$-valued random variable $X$ such that the random configuration $T^{-X}\Gamma$ has law $\mu^*$. An extra head scheme is called *nonrandomized* if it is almost surely equal to a deterministic function of the configuration.

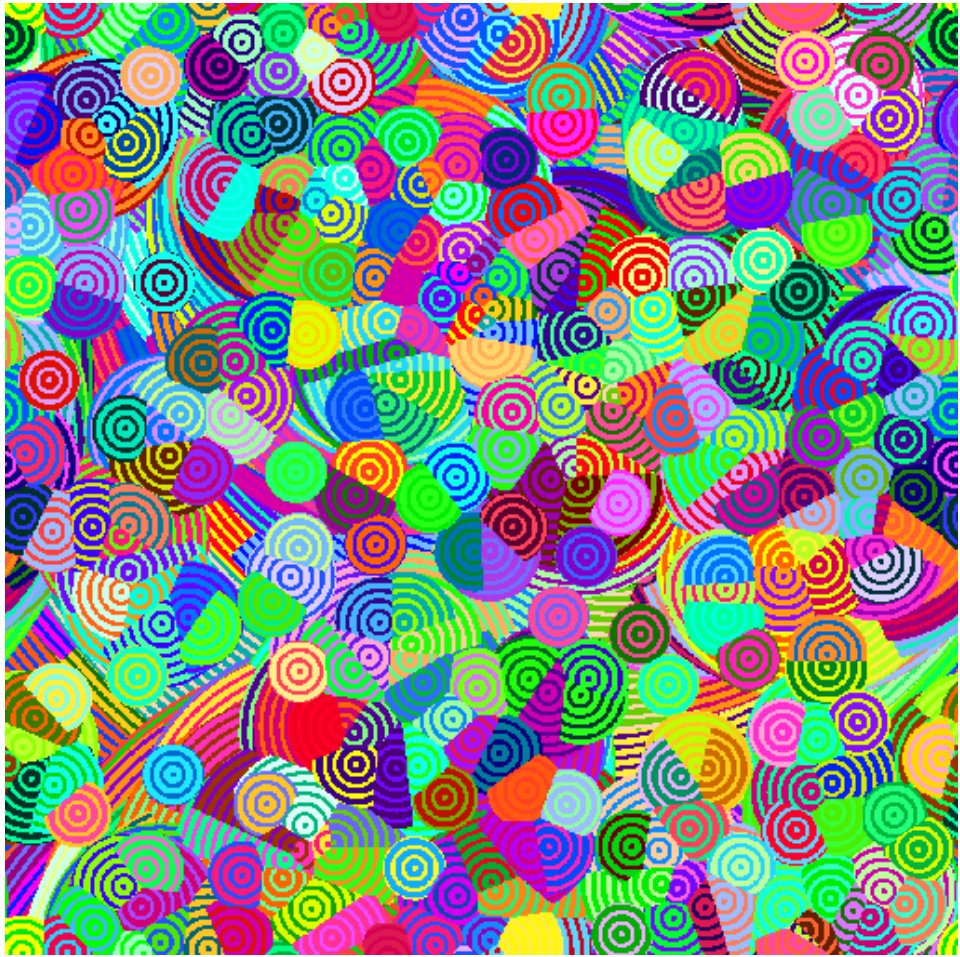

Fig. 1. *A balanced allocation rule applied to a two-dimensional Poisson process (here on a torus). The points of the process are the centers of the concentric circles. Each center is allocated exactly one unit of area, indicated by concentric anulli in two colors. (If you are looking at a greyscale image, color versions are available at www.math.ubc.ca/~holroyd/stable.html and via arXiv:math.PR/0306402.)*

THEOREM 3. *Let $d \geq 1$ and let $\mu$ be an ergodic translation-invariant measure on $\{0,1\}^{\mathbb{Z}^d}$.*

(i) *For all $d, \mu$, there exists an extra head scheme.*

(ii) *For all $d$, there exists a nonrandomized extra head scheme if and only if the marginal probability $p$ is the reciprocal of an integer.*

(iii) *For all $d, \mu$ and any shift coupling of $\Gamma, \Gamma^*$ where $\Gamma = T^X \Gamma^*$, the translation $X$ cannot be a function of $\Gamma^*$.*

Thorisson [13] proved Theorem 3(i). The "if" part of (ii) follows from [8], where appropriate nonrandomized extra head schemes are constructed. We shall present a construction which gives extra head schemes for all $d, \mu$, and also extends to arbitrary countable groups in place of $\mathbb{Z}^d$. When $p$ is rational our construction will yield an extra head scheme which is a deterministic function of $\Gamma$ and an independent roll of a $u$-sided die, where $u$ is the *numerator* of $p$ expressed in its lowest terms.

Consider now the special case when $\mu$ is product measure with parameter $p \in (0, 1)$. It is natural to ask how large the random variable $\|X\|$ must be when $X$ is an extra head scheme (where $\|\cdot\|$ is the Euclidean norm, say). This was essentially answered in dimensions $d = 1, 2$ by Liggett [8] and Holroyd and Liggett [5].

THEOREM 4 ([8], $d = 1$ and [5], $d \geq 2$).  *Let $\mu$ be product measure with parameter $p$ on $\mathbb{Z}^d$.*

(i) *For all $d$, there exists an extra head scheme $X$ satisfying*
$$\mathbf{P}(\|X\| > r) \leq c r^{-d/2},$$
*where $c = c(d, p) < \infty$.*

(ii) *For $d = 1, 2$, any extra head scheme satisfies*
$$\mathbf{E}\|X\|^{d/2} = \infty.$$

It was also shown in [5] that for all $d \geq 1$, any extra head scheme must involve the *examination* of sites at distance at least $Z$ from $O$, where $\mathbf{P}(Z > r) \geq c'(d, p) r^{-d/2}$. In the light of the above results, one might guess that any extra head scheme must satisfy $\mathbf{E}\|X\|^{d/2} = \infty$ for $d \geq 3$ also. In fact, this is very far from the truth.

THEOREM 5.  *Let $\mu$ be product measure with parameter $p$ on $\mathbb{Z}^d$. If $d \geq 3$, then there exists an extra head scheme satisfying*
$$\mathbf{E} \exp(C\|X\|^d) < \infty$$
*for some $C = C(d, p) > 0$.*

(An analogous result also applies to continuum extra head schemes for the Poisson process in $d \geq 3$.) The above result is the best possible up to the value of $C$. Indeed, if $X$ is an extra head scheme, then $\|X\|$ must be at least as large as the distance to the closest occupied site to the origin, so $\mathbf{P}(\|X\| > r) \geq \exp(-C'r^d)$ for some $C' = C'(d, p) > 0$. The proof of Theorem 5 relies on a result of Talagrand [11] on transportation cost.

Consider now the case when $d = 1$ and $\mu$ is an ergodic translation-invariant measure on $\{0, 1\}^{\mathbb{Z}}$. The following natural measure-theoretic construction of an extra head scheme is due to Thorisson [12, 13], and is also presented in [8]. For two measures $\alpha, \beta$ on $\{0,1\}^{\mathbb{Z}}$, define $\alpha \wedge \beta$ to be the measure whose Radon–Nikodym derivative with respect to $\alpha + \beta$ is the minimum of the Radon–Nikodym derivatives of $\alpha$ and $\beta$ with respect to $\alpha + \beta$. Define measures $\alpha_n, \beta_n, \chi_n$ on $\{0,1\}^{\mathbb{Z}}$ for $n \geq 0$ as follows:

$$\alpha_0 = \mu, \qquad \beta_0 = \mu^*$$

and for $n \geq 0$:

$$\chi_n = \alpha_n \wedge T^n \beta_n, \qquad \alpha_{n+1} = \alpha_n - \chi_n, \qquad \beta_{n+1} = \beta_n - T^{-n} \chi_n.$$

Let $X^{\text{meas}}$ be such that

$$\mathbf{P}(X^{\text{meas}} = n, \Gamma \in A) = \chi_n(A).$$

It follows from results in [12, 13] that $X^{\text{meas}} < \infty$ and that $X^{\text{meas}}$ is an extra head scheme. However, the above description gives little clue about how to explicitly construct $X^{\text{meas}}$ from the configuration $\Gamma$.

In contrast, the extra head schemes described in [8] for $\mathbb{Z}$ involve an explicit construction of $X$ from $\Gamma$, and this construction enabled computation of tail behavior. Liggett [8] commented that such solutions were "completely different" from $X^{\text{meas}}$ above. In fact, it turns out that they are identical when $p$ is the reciprocal of an integer. Moreover, we can give a simple explicit construction of $X^{\text{meas}}$ for general $p$.

Let $\Gamma$ have law $\mu$, and let $U$ be a Uniform$(0,1)$ random variable, independent of $\Gamma$. Define $X^{\text{walk}}$ by

$$X^{\text{walk}} = \min\left\{n \geq 0 : \sum_{i=0}^{n}(1 - p^{-1}\Gamma(i)) < U\right\}.$$

(See Figure 2.)

PROPOSITION 6. $X^{\text{meas}}$ and $X^{\text{walk}}$ are extra head schemes, and the joint laws of $(X^{\text{meas}}, \Gamma)$ and $(X^{\text{walk}}, \Gamma)$ are identical.

It is easy to check that $X^{\text{walk}}$ is the same as the extra head scheme constructed by Liggett [8] when $p$ is the reciprocal of an integer.

Our main tool will be a bijective correspondence between extra head schemes and *balanced transport rules* (to be defined later). In the special case of nonrandomized extra head schemes, the correspondence becomes simpler, and can be expressed instead in terms of *balanced allocation rules*. We describe this case below.

Let $\mu$ be a translation-invariant ergodic measure on $\{0,1\}^{\mathbb{Z}^d}$, and suppose that the marginal probability $p$ is the reciprocal of an integer. A (*discrete*) *balanced allocation rule* for $\mu$ is a measurable map $\Phi$ which assigns to $\mu$-almost-every configuration $\gamma$ and every site $x$ a site $\Phi_\gamma(x)$, such that the following properties hold. First, we have $|(\Phi_\gamma)^{-1}(y)| = p^{-1}\gamma(y)$ for $\mu$-almost-all $\gamma$ and all $y$; that is, almost surely the range of $\Phi_\Gamma$ is the set of occupied sites, and each occupied site has exactly $p^{-1}$ pre-images. Second, $\Phi$ is translation-invariant in the sense that if $\Phi_\gamma(x) = y$, then $\Phi_{T^z\gamma}(T^z x) = T^z y$.

PROPOSITION 7. *Let $\Gamma$ have law $\mu$, and suppose $p$ is the reciprocal of an integer. If $\Phi$ is a balanced allocation rule for $\mu$, then the random variable $X$ given by*

$$X = \Phi_\Gamma(0) \tag{1}$$

*is a nonrandomized extra head scheme for $\mu$. Conversely, if $X$ is a nonrandomized extra head scheme, then there exists a $\mu$-almost-everywhere unique balanced allocation rule $\Phi$ satisfying* (1).

Suppose that $p = \frac{1}{2}$ and consider the natural special case of a nonrandomized extra head scheme $X$ such that $X = 0$ whenever $\Gamma(0) = 1$. We call such an $X$ *lazy*. This corresponds via Proposition 7 to a balanced allocation rule $\Phi$ in which for every occupied site $x$ we have $\Phi_\Gamma(x) = x$ almost surely. Such a $\Phi$ amounts to an translation-invariant *matching rule* of occupied sites to unoccupied sites, in which unoccupied site $x$ is matched to occupied site $\Phi_\Gamma(x)$. Then $X$ equals the origin if it is occupied, or the partner of the origin otherwise.

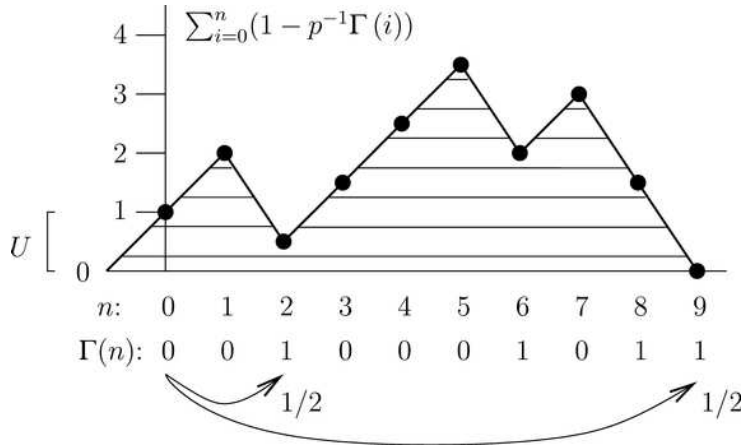

FIG. 2. *An illustration of the construction of $X^{\mathrm{walk}}$. The walk $\sum_{i=0}^n (1 - p^{-1}\Gamma(i))$ is plotted as a function of $n$. In this example $p = \frac{2}{5}$, so the walk takes an up-step of $1$ for an unoccupied site and a down-step of $\frac{3}{2}$ for an occupied site. Conditional on this configuration, $X^{\mathrm{walk}}$ takes the values $2, 9$ each with probability $\frac{1}{2}$.*

We shall use Proposition 7 and its generalizations to deduce results about extra head schemes from results about allocations. The reverse implication is also potentially useful. As an illustration, we note that the following are consequences of Theorem 4(ii) combined with our results.

COROLLARY 8. *Let $\mu$ be product measure on $\mathbb{Z}^d$ with parameter $p$ the reciprocal of an integer. If $d = 1, 2$, then any balanced allocation rule $\Phi$ for $\mu$ satisfies*

$$\mathbf{E}\|\Phi_\Gamma(0)\|^{d/2} = \infty.$$

We shall also state a continuum analogue of Corollary 8.

COROLLARY 9. *Let $\Pi, \Pi'$ be two independent Poisson processes of unit intensity in $\mathbb{R}^d$ and consider any translation-invariant random perfect matching between the points of $\Pi$ and the points of $\Pi'$. If $d = 1, 2$, then the total distance $L$ from points in $[0,1)^d$ to the points they are matched to satisfies*

$$\mathbf{E}L^{d/2} = \infty.$$

Consider the extra head scheme $X^{\text{walk}}$ in Proposition 6 when $d = 1$ and $p = \frac{1}{2}$. Note $X^{\text{walk}}$ is lazy, therefore it corresponds to a matching rule. It is easy to see that the matching rule has the following simple description. Wherever the sequence $\ldots, \Gamma(-1), \Gamma(0), \Gamma(1), \ldots$ has an adjacent pair of the form $(\Gamma(i), \Gamma(i+1)) = (0, 1)$, match them to each other. Then remove all such pairs from the sequence and repeat indefinitely. This matching was used earlier by Meshalkin [10] in the context of finitary isomorphisms.

When $d = 1$, one might guess that $X^{\text{meas}}$ is optimal in the sense that any other nonnegative extra head scheme stochastically dominates it; Srinivasa Varadhan asked whether this was the case (personal communication). The answer is no. For a counterexample, let $\mu$ be product measure with parameter $\frac{1}{2}$. Wherever the configuration contains a sequence of the form $(\Gamma(i), \ldots, \Gamma(i+3)) = (0, 0, 1, 1)$, the allocation rule (Meshalkin matching) corresponding to $X^{\text{walk}} = X^{\text{meas}}$ above has $\Phi_\Gamma(i) = i + 3$ and $\Phi_\Gamma(i+1) = i + 2$. Consider modifying the matching rule so that instead $\Phi_\Gamma(i) = i + 2$ and $\Phi_\Gamma(i+1) = i + 3$ in this situation. By Proposition 7 this results in an extra head scheme $X'$ satisfying $\mathbf{P}(X' \leq 2) > \mathbf{P}(X^{\text{meas}} \leq 2)$, so $X^{\text{meas}}$ was not stochastically optimal. On the other hand, one may similarly show (by induction) that no nonnegative extra head scheme can be strictly stochastically dominated by $X^{\text{walk}}$. Hence there is no stochastically optimal extra head scheme.

The article is organized as follows. In Sections 2 and 3 we establish correspondences of extra head schemes with transports and allocations, and prove

Proposition 7 and Corollaries 8 and 9. In Sections 4 and 5 we construct allocations and transports, and prove Theorem 1 and Theorem 3(i) and (ii). In Section 6 we prove Proposition 2 and Theorem 3(iii) regarding shift coupling in the reverse direction. In Section 7 we prove Proposition 6 about one-dimensional constructions, and in Section 8 we prove the tail estimate Theorem 5.

**2. Discrete equivalence.** In this section we state and prove an equivalence between discrete extra head schemes and balanced transport rules, of which Proposition 7 is a special case.

Let $G$ be an infinite countable group with identity $\mathbf{i}$, and let $\mu$ be a measure on the product $\sigma$-algebra of $\{0,1\}^G$. Elements of $G$ are called *sites* and elements of $\{0,1\}^G$ are called *configurations*. A site $g$ acts on other sites $x$ via left multiplication $g: x \mapsto gx$, and hence on configurations via $(g\gamma)(x) = \gamma(g^{-1}x)$, on measurable functions $f:\{0,1\}^G \to \mathbb{R}$ via $(gf)(\gamma) = f(g^{-1}\gamma)$, on events $A \subseteq \{0,1\}^G$ via $gA = \{g\gamma : \gamma \in A\}$ (whence $\mathbf{1}[gA] = g\mathbf{1}[A]$ where $\mathbf{1}[A]$ denotes the indicator of $A$), and on measures via $(g\mu)(f) = \mu(g^{-1}f)$. We suppose that $\mu$ is invariant and ergodic under the action of $G$. We write $p$ for the marginal probability

$$p = \mu(\Gamma(\mathbf{i}) = 1).$$

We assume that $0 < p < 1$, and we write $\mu^*$ for the conditional law of $\Gamma$ given $\Gamma(\mathbf{i}) = 1$:

$$\mu^*(\cdot) = \mu(\Gamma \in \cdot \,|\, \Gamma(\mathbf{i}) = 1).$$

Let $X$ be a discrete $G$-valued random variable on some joint probability space with $\Gamma$, with probability measure $\mathbf{P}$ and expectation operator $\mathbf{E}$. We call $X$ a (*discrete*) *extra head scheme* for $\mu$ if $X^{-1}\Gamma$ has law $\mu^*$ under $\mathbf{P}$.

A (*discrete*) *transport rule* for $\mu$ is a measurable function $\Theta$ which assigns to $\mu$-almost-every configuration $\gamma$ and every pair of sites $x, y$ a nonnegative real number $\Theta_\gamma(x, y)$, with properties (2), (3) as follows. We think of $\Theta_\gamma(x, y)$ as the amount of mass transported from $x$ to $y$ when the configuration is $\gamma$, and we write

$$\Theta_\gamma(A, B) = \sum_{x \in A, y \in B} \Theta_\gamma(x, y).$$

We require the following properties. First,

(2) $$\Theta_\Gamma(x, G) = 1$$

for $\mu$-almost-all $\Gamma$ and all $y$ (i.e., each site sends out exactly one unit in total). Second, $\Theta$ is $G$-invariant in the sense that

(3) $$\Theta_{g\gamma}(gx, gy) = \Theta_\gamma(x, y)$$

for all $\gamma$ and all $x, y, g \in G$.

We call a transport rule $\Theta$ *balanced* if it satisfies in addition

$$\Theta_\Gamma(G, y) = p^{-1}\Gamma(y), \tag{4}$$

for $\mu$-almost-all $\Gamma$ and all $x$, $y$ (i.e., unoccupied sites receive nothing, all occupied sites receive equal mass, which must then necessarily be $p^{-1}$).

We are now ready to state the equivalence result. Fix $\mu$, let $\Theta$ be a transport rule, let $X$ be a $G$-valued random variable and suppose that

$$\Theta_\Gamma(\mathbf{i}, x) = \mathbf{P}(X = x|\Gamma) \tag{5}$$

for $\mu$-almost-all $\Gamma$ and all $x$ (i.e., conditional on the configuration, the identity distributes one unit of mass according to the conditional distribution of $X$). Note that by summing over $x$ and using (3), (5) implies (2). For any $X$, (5) determines $\Theta$ uniquely up to a **P**-null event, and conversely for any $\Theta$, (5) uniquely determines the joint law of $X$, $\Gamma$.

THEOREM 10. *Suppose that $X$ and $\Theta$ are related by* (5). *Then $X$ is an extra head scheme if and only if $\Theta$ is balanced.*

PROOF OF PROPOSITION 7. This is an immediate special case of Theorem 10, where $G$ is $\mathbb{Z}^d$ under addition, and we identify a balanced allocation rule $\Phi$ with the balanced transport rule given by $\Theta_\gamma(x, y) = \mathbf{1}[\Phi_\gamma(x) = y]$. □

PROOF OF COROLLARY 8. Immediate from Proposition 7 and Theorem 4. □

We shall make use of the following lemma.

LEMMA 11 (Mass transport principle). *Let $m: G \times G \to [0, \infty]$ be such that $m(gx, gy) = m(x, y)$ for all $x$, $y$, $g$. Then*

$$\sum_{y \in G} m(x, y) = \sum_{y \in G} m(y, x).$$

For a proof see [1] or [3].

The proof of Theorem 10 is based on the following lemma. Let $J$ be the total mass received by the identity:

$$J = J(\Gamma) = \Theta_\Gamma(G, \mathbf{i}).$$

LEMMA 12. *Suppose $X$ and $\Theta$ are related by (5). For any nonnegative measurable function $f$ on $\{0,1\}^G$, we have*
$$\mathbf{E}(f(X^{-1}\Gamma)) = \mathbf{E}(J(\Gamma)f(\Gamma)).$$

[In the above, $J(\Gamma)f(\Gamma)$ denotes ordinary multiplication.]

PROOF OF LEMMA 12. The following device will be useful. Enlarging the probability space if necessary, we may assume that $X$ is a deterministic function of $\Gamma$ and an independent Uniform$(0,1)$ random variable $U$; thus, $X = \xi(\Gamma, U)$. ($U$ represents any "additional randomization" in the choice of $X$; see [5] for a more detailed explanation.)

We have the chain of equalities
$$\mathbf{E}(f(X^{-1}\Gamma)) = \int d\mu(\gamma) \int_0^1 du\, f(\xi(\gamma, u)^{-1}\gamma)$$
$$= \int d\mu(\gamma) \int_0^1 du \sum_{x \in G} \mathbf{1}[\xi(\gamma, u) = x] f(x^{-1}\gamma)$$
(6)
$$= \int d\mu(\gamma) \sum_{x \in G} \Theta_\gamma(\mathbf{i}, x) f(x^{-1}\gamma)$$
(7)
$$= \int d\mu(\gamma)\, \Theta_\gamma(G, \mathbf{i}) f(\gamma)$$
$$= \mathbf{E}(J(\Gamma)f(\Gamma)).$$

In (6) we have used (5), and in (7) we have used Lemma 11 with $m(x,y) = \mathbf{E}\Theta_\Gamma(x,y)f(y^{-1}\Gamma)$. □

PROOF OF THEOREM 10. Suppose that $\Theta$ is a balanced transport rule. For any nonnegative measurable $f$, by Lemma 12 and (4) we have
$$\mathbf{E}(f(X^{-1}\Gamma)) = \mathbf{E}(J(\Gamma)f(\Gamma)) = p^{-1}\mathbf{E}(\Gamma(\mathbf{i})f(\Gamma))$$
$$= \mathbf{E}(f(\Gamma)|\Gamma(\mathbf{i}) = 1) = \mu^*(f).$$
So $X^{-1}\Gamma$ has law $\mu^*$, thus $X$ is an extra head scheme.

Conversely, suppose that $X$ is an extra head scheme. We must check that $\Theta$ is balanced. Since $(X^{-1}\Gamma)(\mathbf{i}) = 1$ almost surely, it is immediate from (5) and (3) that every unoccupied site receives zero mass, so it is sufficient to check that every occupied site receives mass $p^{-1}$ almost surely. By (3) it is enough to check this for $\mathbf{i}$, so we must check that under $\mu^*$ we have $J = p^{-1}$ almost surely.

Since $X$ is an extra head scheme, for any $f$ we have $\mathbf{E}(f(X^{-1}\Gamma)) = \mu^*(f)$. Note also that $\mathbf{E}(Jf) = p\mathbf{E}(Jf|\Gamma(\mathbf{i}) = 1) + (1-p)\mathbf{E}(Jf|\Gamma(\mathbf{i}) = 0) = p\mu^*(Jf)$ [since $J = 0$ on $\{\Gamma(\mathbf{i}) = 0\}$]. Thus Lemma 12 yields
$$\mu^*(f) = p\mu^*(Jf).$$

Applying this first with $f \equiv 1$ and then with $f = J$ shows that under $\mu^*$, the random variable $J$ has mean $p^{-1}$ and variance 0; hence, $\mu^*$-almost-surely we have $J = p^{-1}$. □

**3. Continuum equivalence.** The equivalence between extra head schemes and balanced transport rules in Theorem 10 has an analogue in the continuum setting, which we shall state (without proof) at the end of this section. Since the full continuum result is somewhat technical and is not required for any of our main results, we shall instead prove the special case involving nonrandomized extra head schemes and allocations (the analogue of Theorem 7).

Let $\Pi$ be a translation-invariant ergodic simple point process of intensity 1 on $\mathbb{R}^d$, with law $\Lambda$. Elements of $\mathbb{R}^d$ are called *sites*. Integer-valued Borel measures on $\mathbb{R}^d$ are called *configurations*. Let $\mathcal{L}$ denote Lebesgue measure on $\mathbb{R}^d$. For any $z \in \mathbb{R}^d$, we define the translation $T^z$, which acts on sites via $T^z x = x + z$, on functions $h : \mathbb{R}^d \to \mathbb{R}$ via $(T^z h)x = h(T^{-z}x)$ and on configurations via $(T^z \pi)(h) = \pi(T^{-z}h)$.

Let $\Pi^*$ be the Palm version of $\Pi$, with law $\Lambda^*$. The following is a standard property of the Palm process. For any bounded measurable function $f$ on configurations and any Borel set $B \subseteq \mathbb{R}^d$, we have

(8) $$\mathbf{E} \int_B f(T^{-s}\Pi) \Pi(ds) = \mathcal{L}(B) \mathbf{E} f(\Pi^*).$$

Note that the integral on the left-hand side can be written as $\sum_{s \in [\Pi] \cap B} f(T^{-s}\Pi)$. See, for example, [7] for details.

A (*continuum*) *allocation rule* for $\Pi$ is a measurable function $\Psi$ which assigns to $\Lambda$-almost-every configuration $\pi$ and every site $x$ a site $\Psi_\pi(x)$, and which is translation-invariant in the sense that if $\Psi_\pi(x) = y$, then $\Psi_{T^z\pi}(T^z x) = T^z y$. (It is important that we require the preceding statement to hold for *all* configurations $\pi$; in particular, it is thus understood that $\Psi_{T^z\pi}$ is defined if and only if $\Psi_\pi$ is.) Let $\mathcal{L}$ denote Lebesgue measure on $\mathbb{R}^d$. An allocation rule $\Psi$ is called *balanced* if $\Lambda$-almost-surely for each $s \in [\Pi]$ we have $\mathcal{L}(\Psi_\Pi^{-1}(s)) = 1$, while $\mathcal{L}(\Psi_\Pi^{-1}(\mathbb{R}^d \setminus [\Pi])) = 0$.

THEOREM 13. *Let $\Psi$ be an allocation rule for $\Pi$. The random variable $Y = \Psi_\Pi(0)$ is a nonrandomized extra head scheme for $\Pi$ if and only if $\Psi$ is balanced.*

We shall prove Theorem 13 via Lemma 14. Let $\Psi$ be an allocation rule. For $z \in \mathbb{Z}^d$, define the unit cube $Q_z = z + [0,1)^d \subseteq \mathbb{R}^d$. For $s \in \mathbb{R}^d$, write $J_\Pi(s) = \mathcal{L}(\Psi_\Pi^{-1}(s))$ and $\Pi_s = T^{-\Psi_\Pi(s)}\Pi$.

LEMMA 14. *For any $z \in \mathbb{Z}^d$ and any nonnegative measurable $f$, we have*

$$\mathbf{E}f(\Pi_0) = \mathbf{E} \int_{Q_z} J_\Pi(s) f(T^{-s}\Pi) \Pi(ds).$$

PROOF. The translation-invariance of $\Lambda$ and $\Psi$ implies that $\Pi_x$ has the same law for each $x \in \mathbb{R}^d$. Indeed, write $\Pi'$ for $T^{-x}\Pi$. Then $\Psi_\Pi(x) = x + \Psi_{\Pi'}(0)$, so that $T^{-\Psi_\Pi(x)}\Pi = T^{-\Psi_{\Pi'}(0)}(\Pi') = \Pi'_0$, which has the law of $\Pi_0$. Therefore, $\mathbf{E}f(\Pi_0) = \mathbf{E}f(\Pi_x)$ for any $f$. Fix $f$ and $x$, and define

$$m(z,w) = \mathbf{E} \int_{Q_z} f(\Pi_x) \mathbf{1}[\Psi_\Pi(x) \in Q_w] \mathcal{L}(dx).$$

Applying the mass transport principle (Lemma 11) yields

$$\sum_{w \in \mathbb{Z}^d} m(z,w) = \sum_{w \in \mathbb{Z}^d} m(w,z).$$

The left-hand side equals $\mathbf{E}f(\Pi_0)$, and the right-hand side equals

$$\mathbf{E} \int_{\mathbb{R}^d} f(T^{-\Psi_\Pi(x)}\Pi) \mathbf{1}[\Psi_\Pi(x) \in Q_z] \mathcal{L}(dx) = \mathbf{E} \int_{Q_z} J_\Pi(s) f(T^{-s}\Pi) \Pi(ds). \quad \square$$

PROOF OF THEOREM 13. If $\Psi$ is balanced, then Lemma 14 immediately gives that $\Psi_\Pi(0)$ is an extra head scheme. For the converse, apply the lemma to $f \equiv 1$ and $f(\pi) = J_\pi(0)$. $\quad \square$

The following is the continuum analogue of Corollary 8.

COROLLARY 15. *Let $\Pi$ be a Poisson process of unit intensity on $\mathbb{R}^d$. If $d = 1, 2$, then any balanced allocation rule $\Psi$ for $\Pi$ satisfies*

$$\mathbf{E}\|\Psi_\Pi(0)\|^{d/2} = \infty.$$

PROOF. One possible proof is to deduce the result from Theorem 13 together with Theorem 2(B) of [5], which is the continuum analogue of Theorem 4(ii). Alternatively, we may proceed via discrete transports as follows.

Denote the unit cube $Q_z = z + [0,1)^d$. Let $\Psi$ be a balanced allocation rule for $\Pi$, and define a discrete configuration $\Gamma$ by

$$\Gamma(z) = 1 \wedge \Pi(Q_z),$$

so that the law of $\Gamma$ is product measure with parameter $1 - e^{-1}$ on $\mathbb{Z}^d$. Now define $\Theta$ by

$$\Theta_\Gamma(x,y) = \mathbf{E}(\mathcal{L}[\Psi_\Pi^{-1}(Q_y) \cap Q_x] | \Gamma).$$

It is elementary to check that $\Theta$ is a balanced transport rule for $\Gamma$, so by Theorem 10 there is an associated extra head scheme $X$. It is furthermore easy to check that $\mathbf{E}\|\Psi_\Pi(0)\|^{d/2} < \infty$ implies $\mathbf{E}\|X\|^{d/2} < \infty$, so the result follows from Theorem 4(ii). □

PROOF OF COROLLARY 9. The required statement may be formulated as follows. Let $M$ be a simple point process on $\mathbb{R}^d \times \mathbb{R}^d$, invariant under the diagonal action of translations of $\mathbb{R}^d$. We write $M(A,B) = M(A \times B)$, and suppose that the marginals given by $\Pi(\cdot) = M(\mathbb{R}^d, \cdot)$ and $\Pi'(\cdot) = M(\cdot, \mathbb{R}^d)$ are two independent Poisson processes of unit intensity on $\mathbb{R}^d$. [If $M$ has an atom at $(x,y)$, it means that the point $x$ of $\Pi'$ is matched to the point $y$ of $\Pi$.] It is sufficient to prove that for $d = 1, 2$, any such $M$ must satisfy

$$(9) \qquad \iint \|x-y\|^{d/2} \mathbf{1}[x \in Q_0] M(dx, dy) = \infty.$$

As in the preceding proof, we define

$$\Gamma(z) = 1 \wedge \Pi(Q_z),$$

and

$$\Theta_\Gamma(x,y) = \mathbf{E}(M(Q_x, Q_y)|\Gamma).$$

It is easy to check that the law of $\Gamma$ is product measure on $\mathbb{Z}^d$, and that $\Theta$ is a balanced transport for $\Gamma$. Equation (9) may then be deduced from Theorem 10 and Theorem 4(ii). □

Finally in this section we shall state without proof the full continuum analogue of Theorem 10. A *transport* is a nonnegative $\sigma$-finite Borel measure $\omega$ on $\mathbb{R}^d \times \mathbb{R}^d$. We write $\omega(A,B) = \omega(A \times B)$, and think of this as the mass sent from $A$ to $B$. The *marginals* of $\omega$ are the measures $\omega(\cdot, \mathbb{R}^d), \omega(\mathbb{R}^d, \cdot)$ on $\mathbb{R}^d$. Let $\Pi$ be a translation-invariant, ergodic simple point process on $\mathbb{R}^d$ with law $\Lambda$. A (*continuum*) *transport rule* for $\Pi$ is a measurable map $\Omega$ which assigns to $\Lambda$-almost-every configuration $\pi$ a transport $\Omega_\pi$, with the following properties. The first marginal $\Omega_\Pi(\cdot, \mathbb{R}^d)$ is Lebesgue measure $\Lambda$-almost-surely, and $\Omega$ is invariant in the sense that $\Omega_{T^z\pi}(T^zA, T^zB) = \Omega_\pi(A,B)$ for all $\pi, z, A, B$. A transport rule $\Omega$ is *balanced* if the second marginal satisfies $\Lambda$-almost-surely $\Omega_\Pi(\mathbb{R}^d, A) = \Pi(A)$ for all $A \subseteq \mathbb{R}^d$.

Let $Y$ be an $\mathbb{R}^d$-valued random variable and let $\Omega$ be a transport rule, and suppose $\mathbf{P}$ admits conditional probabilities such that

$$(10) \qquad \frac{d[\Omega_\Pi(\cdot, A)]}{d\mathcal{L}(\cdot)}(0) = \mathbf{P}(Y \in A|\Pi).$$

Here a specific version of the Radon–Nikodym derivative must be used, to ensure that it is defined everywhere and translation-invariant. By the Lebesgue differentiation theorem (see [9], Theorem 2.1.2), the upper density $\limsup_{r \to 0} \nu(B(x,r))/\mathcal{L}(B(x,r))$ is a suitable version of the Radon–Nikodym derivative $d\nu/d\mathcal{L}$.

THEOREM 16. *Suppose $Y$ and $\Omega$ are related as in* (10). *Then $Y$ is an extra head scheme if and only if $\Omega$ is balanced.*

We omit the proof of Theorem 16, which proceeds along the same lines as that of Theorem 10. The proof involves no new ideas, but more technical notation.

**4. Discrete allocations and transports.** Let $\mu$ be an ergodic $G$-invariant measure on $\{0,1\}^G$. In this section we shall prove the following.

THEOREM 17. *For any $G$, $\mu$, there exists a balanced transport rule.*

THEOREM 18. *For any $G$, there exists an integer-valued balanced transport rule if and only if $p$ is the reciprocal of an integer.*

THEOREM 19. *For any $G, \mu$, there exists an extra head scheme.*

THEOREM 20. *For any $G$, there exists a nonrandomized extra head scheme if and only if $p$ is the reciprocal of an integer.*

PROOF OF THEOREM 18, "ONLY IF" PART. In an integer-valued transport rule, the unit of mass sent out by a site all goes to a single site, while in a balanced transport rule, occupied sites receive mass $p^{-1}$. Hence $p^{-1}$ must be an integer. □

PROOF OF THEOREM 3(i), (ii) AND THEOREMS 19, 20. Theorems 19 and 20 follow immediately from Theorems 17 and 18 together with Theorem 10. [A nonrandomized extra head scheme corresponds via (5) to an integer-valued balanced transport rule.] Theorem 3(i), (ii) are Theorems 19 and 20 specialized to $\mathbb{Z}^d$. □

PROOF OF THEOREM 17. We construct the required transport rule by a kind of invariant greedy algorithm. Order the elements of $G$ as $G = \{g_0, g_1, \ldots\}$, and fix a configuration $\gamma$. Informally, each site starts with mass 1 to distribute, while a site $y$ has the capacity to receive mass $p^{-1}\gamma(y)$. At time $n$, every site $x$ sends as much mass as possible to site $g_n x$. Formally, inductively define $\theta^n(x,y) = \theta^n_\gamma(x,y)$ for $n = 0, 1, \ldots$ as follows. For all sites $x, y$,

$$\theta^0(x,y) = 0,$$

and for $n \geq 0$,

$$\theta^{n+1}(x,y) = \theta^n(x,y) + \delta^n(x,y),$$

where

$$\delta^n(x,y) = \mathbf{1}[g_n x = y] \min\{1 - \theta^n(x,G), p^{-1}\gamma(y) - \theta^n(G,y)\}.$$

Finally, put $\Theta_\gamma(x,y) = \lim_{n\to\infty} \theta_{\gamma^n}(x,y)$.

Clearly, $\Theta$ is $G$-invariant; we claim that it is a balanced transport rule. By the construction, we have for all $x$

$$\Theta_\Gamma(x,G) \leq 1 \quad \text{and} \quad \Theta_\Gamma(G,x) \leq p^{-1}\Gamma(x).$$

We call a site $x$ *unexhausted* if the former inequality is strict, and we call $x$ *unsated* if the latter inequality is strict. We must show that almost surely there are no unexhausted sites and no unsated sites. First, note that unexhausted sites and unsated sites cannot exist simultaneously for the same $\gamma$. For suppose that $x$ is unexhausted and $y$ is unsated. Then considering $\delta^n(x,y)$ where $n$ is such that $g_n x = y$ shows that either $\theta^{n+1}(x,y) = 1$ or $\theta^{n+1}(x,y) = p^{-1}\gamma(y)$, a contradiction. Also, by ergodicity, the existence of unexhausted sites and the existence of unsated sites are both zero–one events. Hence it remains only to rule out the possibility that almost surely one occurs without the other. The mass transport principle (Lemma 11) applied to $m(x,y) = \mathbf{E}\Theta_\Gamma(x,y)$ yields

$$\mathbf{E}\Theta_\Gamma(0,G) = \mathbf{E}\Theta_\Gamma(G,0),$$

but the left-hand side is less that 1 if and only if there exist unexhausted sites, and the right-hand side is less that 1 if and only if there exist unsated sites. □

REMARK. In the case when $G = \mathbb{Z}$ under addition, the above construction also gives a balanced transport rule if we set $g_n = n$ for all $n \geq 0$, even though $g_0, g_1, \ldots$ no longer exhausts $G$. (This will be relevant in the proof of Proposition 6.) The above proof goes through, except for the argument that unexhausted sites and unsated sites cannot exist simultaneously, which must be modified as follows. By the previous argument, for $x \leq y$ it is impossible that $x$ is unexhausted and $y$ is unsated. Hence if with positive probability both unexhausted and unsated sites existed, then by the invariance of the construction, the random variable $\max\{x : x \text{ is unsated}\}$ would take all integer values with equal positive probabilities, which is impossible.

PROOF OF THEOREM 18, "IF" PART. Consider the construction of $\Theta$ in the proof of Theorem 17 above. If $p^{-1}$ is an integer, then each $\theta^n$ is integer-valued, so the same applies to $\Theta$. □

If $p = u/v$ where $u,v$ are integers, the same argument shows that $u^{-1}\Theta$ is integer-valued, and the corresponding extra head scheme can consequently be written as a deterministic function of $\Gamma$ and an independent roll of a $u$-sided die, as remarked in the Introduction.

Note also that if the ordering of $G$ satisfies $g_0 = \mathbf{i}$, then the resulting extra head scheme is lazy.

**5. Continuum allocations.** Let $\Pi$ be a translation-invariant ergodic simple point process of unit intensity on $\mathbb{R}^d$. Denote the law of $\Pi$ by $\Lambda$.

THEOREM 21. *For any $d, \Lambda$, there exists a balanced continuum allocation rule.*

PROOF OF THEOREM 1. Immediate from Theorems 21 and 13. □

It is natural to try to prove Theorem 21 by some continuous-time version of the invariant greedy algorithm, in which sites of $\mathbb{R}^d$ are ordered by Euclidean norm, say. Although this is an appealing idea, it appears difficult to rigorize directly. Instead, our construction will be based on the stable marriage algorithm of Gale and Shapley [2].

PROOF OF THEOREM 21. In what follows, all distances are Euclidean. Elements of $[\Pi]$ are called $\Pi$-points. Let $L$ be the (random) set of all sites of $\mathbb{R}^d$ which are equidistant from two or more $\Pi$-points. Since $\Pi$ has intensity 1, $[\Pi]$ is countable almost surely; hence, $\mathcal{L}(L) = 0$ almost surely. For convenience we set $\Psi_\Pi(s) = s$ for all $s \in L$.

Consider the following algorithm. For each positive integer $n$, *stage* $n$ consists of two parts as follows.

  (i) Each site $s \notin L$ *applies* to the closest $\Pi$-point to $s$ which has not rejected $s$ at any earlier stage.

  (ii) For each $\Pi$-point $x$, let $A_n(x)$ be the set of sites which applied to $x$ in stage $n$ (i), and define the *rejection radius*

$$r_n(x) = \inf\{r : \mathcal{L}(A_n(x) \cap B(\sigma, r)) \geq 1\},$$

where $B(x, r) = \{s \in \mathbb{R}^d : \|s - x\| < r\}$ is the ball of radius $r$ at $x$, and the infimum of the empty set is taken to be $\infty$. Then $x$ *shortlists* all sites in $A_n(x) \cap B(x, r_n(x))$, and *rejects* all sites in $A_n(x) \setminus B(x, r_n(x))$.

We now describe $\Psi$. Consider a site $s \notin L$. Since any bounded set contains only finitely many $\Pi$-points almost surely, the following is clear. Either $s$ is rejected by every $\Pi$-point (in increasing order of distance from $s$), or, for some $\Pi$-point $x$ and some stage $n$, $s$ is shortlisted by $x$ at stage $n$ and all later stages. In the former case we call $s$ *unclaimed* and put for convenience $\Psi_\Pi(s) = s$; in the latter case we put $\Psi_\Pi(s) = x$.

We claim that $\Psi$ is a balanced allocation rule. Clearly, it satisfies the required measurability and translation-invariance.

Let $S_n(x)$ be the set of sites shortlisted by a $\Pi$-point $x$ at stage $n$. By the construction in (ii) and the intermediate value theorem, we have $\mathcal{L}(S_n(x)) \leq 1$. But by the definition of $\Psi$ above we have $\Psi_\Pi^{-1}(x) = \limsup_{n\to\infty} S_n(x) = \liminf_{n\to\infty} S_n(x)$, so Fatou's lemma implies $\mathcal{L}(\Psi_\Pi^{-1}(x)) \leq 1$. We call a $\Pi$-point $x$ *unsated* if that inequality is strict. Note also that if a $\Pi$-point $x$ ever rejects any sites (at stage $n$, say), then we must have $\mathcal{L}(S_m(x)) = 1$ for all later stages $m \geq n$. Hence an unsated $\Pi$-point never rejected any sites.

We must show that almost surely there are no unsated $\Pi$-points and the set of unclaimed sites is $\mathcal{L}$-null. Unsated $\Pi$-points and unclaimed sites cannot exist simultaneously, since an unclaimed site is rejected by every $\Pi$-point, but an unsated $\Pi$-point never rejects sites. Also, by ergodicity, the existence of unsated $\Pi$-points and of a positive measure of unclaimed sites are both zero–one events, so it remains to rule out the possibility that almost surely one occurs without the other. For $z \in \mathbb{Z}^d$, define the unit cube $Q_z = z + [0,1)^d \subseteq \mathbb{R}^d$. Let

$$m(s,t) = \mathbf{E} \sum_{x \in [\Pi] \cap Q_t} \mathcal{L}(Q_s \cap \Psi_\Pi^{-1}(x)).$$

By the mass transport principle (Lemma 11), we have

$$\mathbf{E} \sum_{x \in [\Pi] \cap Q_0} \mathcal{L}(\Psi_\Pi^{-1}(x)) = \sum_{s \in \mathbb{Z}^d} m(s,0)$$
$$= \sum_{t \in \mathbb{Z}^d} m(0,t) = \mathbf{E}\mathcal{L}(Q_0 \cap \Psi_\Pi^{-1}([\Pi])).$$

Since $\Pi$ has intensity 1, the left-hand side equals 1 if there are no unsated centers, and is strictly less than 1 otherwise. And the right-hand side equals 1 if the set of unclaimed sites is $\mathcal{L}$-null, and is strictly less than 1 otherwise. □

## 6. Reverse extra head schemes.

PROPOSITION 22. *Let $\mu$ be a $G$-invariant ergodic measure on $\{0,1\}^G$, and let $\Gamma$ have law $\mu$. For any discrete extra head scheme $X$ we have almost surely*

$$\mathbf{P}(X = x | X^{-1}\Gamma) \leq p$$

*for all $x \in G$.*

PROPOSITION 23. *Let $\Pi$ be a translation-invariant ergodic point process of unit intensity on $\mathbb{R}^d$. For any continuum extra head scheme $Y$, the conditional law of $Y$ given $T^{-Y}\Pi$ is absolutely continuous with respect to Lebesgue measure, with density bounded above by 1.*

PROOF OF THEOREM 3(iii) AND PROPOSITION 2. Immediate from Propositions 22 and 23. □

PROOF OF PROPOSITION 22. Let $X$ be an extra head scheme for $\Gamma$, and write $\Gamma^* = X^{-1}\Gamma$. Fix $\beta > p$, and define for $x \in G$

$$A_x = \{\gamma^* \in \{0,1\}^G : \mathbf{P}(X = x|\Gamma^*)(\gamma^*) \geq \beta\}.$$

Since $\{\Gamma^* \in A_x, X = x\} \subseteq \{\Gamma \in xA_x\}$, we have

$$\beta\mu^*(A_x) \leq \mathbf{P}(\Gamma^* \in A_x, X = x) \leq \mu(xA_x) = p\mu^*(A_x).$$

Therefore, $\mu^*(A_x) = 0$. Taking a union over rational $\beta > p$ completes the proof.

□

PROOF OF PROPOSITION 23. Let $Y$ be an extra head scheme for $\Pi$, and write $\Pi^* = T^{-Y}\Pi$. It is sufficient to show that for every rational cube $W$ of positive Lebesgue measure, almost surely

$$\mathbf{P}(Y \in W|\Pi^*) \leq \mathcal{L}(W).$$

Fix $\beta > 1$, and define the event

$$A_W = \{\pi^* : \mathbf{P}(Y \in W|\Pi^*)(\pi^*) \geq \beta\mathcal{L}(W)\}.$$

We have

$$\beta\mathcal{L}(W)\Lambda^*(A_W) \leq \mathbf{P}(\Pi^* \in A_W, Y \in W)$$

$$\leq \int \Lambda(d\Pi)\Lambda\left(\bigcup_{y \in [\Pi] \cap W} T^y A_W\right)$$

$$\leq \int \Lambda(d\Pi) \sum_{y \in [\Pi] \cap W} \Lambda(T^y A_W)$$

$$= \mathcal{L}(W)\Lambda^*(A_W).$$

Hence when $\mathcal{L}(W) > 0$, we have $\Lambda^*(A_W) = 0$. □

**7. Measure-theoretic construction.** The following is a variant of the construction of $X^{\mathrm{meas}}$ in the Introduction. Let $\mu$ be a $G$-invariant ergodic measure on $\{0,1\}^G$. Let $G$ be ordered as $G = \{g_0, g_1, \dots\}$. Define measures $\alpha_n, \beta_n, \chi_n$ on $\{0,1\}^G$ as follows:

$$\alpha_0 = \mu, \qquad \beta_0 = \mu^*$$

and for $n \geq 0$:

$$\chi_n = \alpha_n \wedge (g_n\beta_n), \qquad \alpha_{n+1} = \alpha_n - \chi_n, \qquad \beta_{n+1} = \beta_n - g_n^{-1}\chi_n.$$

Let $\Theta_\Gamma$ be the balanced transport rule constructed in the proof of Theorem 17, using the same ordering of $G$ as above. Let $X^{\text{greedy}}$ be the corresponding extra head scheme given by (5) and Theorem 10.

THEOREM 24. *For any $G, \mu$ and any ordering $g_0, g_1, \ldots,$ we have*
$$\mathbf{P}(X^{\text{greedy}} = g_n, \Gamma \in \cdot) = \chi_n(\cdot)$$
*for all $n$.*

PROOF. By construction, the measures $\alpha_n, \beta_n, \chi_n$ are all absolutely continuous with respect to $\mu$. Denote the Radon–Nikodym derivatives
$$a_n = \frac{d\alpha_n}{d\mu}, \qquad b_n = \frac{d\beta_n}{d\mu}, \qquad c_n = \frac{d\chi_n}{d\mu}.$$
We have
$$\alpha_0 = 1, \qquad \beta_0(\gamma) = p^{-1}\gamma(0).$$
And using $G$-invariance of $\mu$,
$$c_n = a_n \wedge (g_n b_n), \qquad a_{n+1} = a_n - c_n, \qquad b_{n+1} = b_n - g_n^{-1} c_n.$$
By induction on $n$, it is easy to verify that
$$a_n(\gamma) = 1 - \theta_\gamma^n(0, G),$$
$$b_n(\gamma) = p^{-1}\gamma(0) - \theta_\gamma^n(G, 0),$$
$$c_n(\gamma) = \delta_\gamma^n(0, g_n) = \Theta_\gamma(0, g_n),$$
where $\theta^n(x, y), \delta^n(x, y)$ are as in the proof of Theorem 17. It follows that for any event $A \subseteq \{0, 1\}^G$,
$$\mathbf{P}(X^{\text{greedy}} = g_n, \Gamma \in A) = \int_A \mu(d\gamma) \Theta_\gamma(0, g_n) = \int_A \mu(d\gamma) c_n(\gamma) = \chi_n(\gamma),$$
as required. □

PROOF OF PROPOSITION 6. Let $g_n = n$ for $n = 0, 1, \ldots$ (note that $g_0, g_1, \ldots$ does not exhaust $G$) and construct $\Theta$, $X^{\text{greedy}}$ and $\chi_n$ as above. As remarked after the proof of Theorem 17, $\Theta$ is a balanced transport rule in this case also, and therefore $X^{\text{greedy}}$ is an extra head scheme. The statement of Proposition 24 above also holds, with the same proof. Therefore, $(X^{\text{greedy}}, \Gamma)$ and $(X^{\text{meas}}, \Gamma)$ have identical joint laws. It remains to check that $(X^{\text{greedy}}, \Gamma)$ and $(X^{\text{walk}}, \Gamma)$ have identical laws. This follows from the fact that for any $x \leq y$,
$$\Theta_\gamma(x, y) = \int_0^1 du\, \mathbf{1}\left[y = \min\left\{z \geq x : \sum_{i=x}^z (1 - p^{-1}\gamma(i)) < u\right\}\right].$$
This is evident from Figure 2. More formally, it may be checked by induction on $y - x$. □

**8. Three-dimensional tail behavior.** In this section we prove Theorem 5.

THEOREM 25. *Let $\mu$ be product measure with parameter $p$ on $\mathbb{Z}^d$. If $d \geq 3$, then there exists a balanced discrete transport rule $\Phi$ satisfying*

$$\mathbf{E}\exp(C\|\Phi_\Gamma(0)\|^d) < \infty$$

*for some $C = C(d,p) > 0$.*

PROOF OF THEOREM 5. Immediate from Theorem 25 and Theorem 10. □

PROOF OF THEOREM 25. It is convenient to work first in a continuum setting, and then transfer to $\mathbb{Z}^d$. A *transport* is a nonnegative $\sigma$-finite measure $\omega$ on $\mathbb{R}^d \times \mathbb{R}^d$. For Borel sets $A, B \subseteq \mathbb{R}^d$ we write $\omega(A, B) = \omega(A \times B)$, and we think of this as the amount of mass transported from $A$ to $B$. By a *random transport* we mean a random element in the space of all transports, this space being equipped with the natural product $\sigma$-algebra. A random transport $\mathcal{T}$ is called *invariant* if $\mathcal{T}(A+z, B+z)$ is equal in law to $\mathcal{T}(A, B)$ for any Borel sets $A$, $B$ and any site $z$. We shall construct an invariant random transport $\mathcal{T}$ whose marginals are Lebesgue measure on $\mathbb{R}^d$ and a Poisson process.

The following is proved in [11]. Let $d \geq 3$. For each positive integer $m$ there exists a random transport $\mathcal{S} = \mathcal{S}_m$ with the following properties. The first marginal $\mathcal{S}(\cdot, \mathbb{R}^d)$ is Lebesgue measure on the cube $[0,1]^d$ almost surely. The second marginal $\mathcal{S}(\mathbb{R}^d, \cdot)$ is equal in law to $m^{-1}\sum_{i=1}^m \delta_i$, where the $\delta_i$ are point masses whose locations are i.i.d. uniform on $[0,1]^d$. Finally, for constants $c, c' < \infty$ depending only on $d$, we have the following "bound on transportation cost":

$$\mathbf{E}\iint \exp(cm\|x-y\|^d)\mathcal{S}(dx, dy) \leq c'.$$

Here $\|\cdot\|$ is the Euclidean norm and $\mathbf{E}$ denotes expectation with respect to the random transport.

We now define a random transport $\mathcal{T}_m$ as follows. Informally, we rescale $\mathcal{S}_m$ to cover a cube of volume $m$, and multiply by $m$ so that the intensity is still 1; then we tile space with identical copies of this transport, with the origin chosen uniformly at random. Formally, let $a$ be uniform on $[0,1]^d$ and independent of $\mathcal{S}_m$, and define $\mathcal{T}_m$ by

$$\mathcal{T}_m(A, B) = m\sum_{z \in \mathbb{Z}^d} \mathcal{S}_m(m^{-1/d}(A+a+z), m^{-1/d}(B+a+z)).$$

It is easy to check the following. $\mathcal{T}_m$ is invariant. $\mathcal{T}_m(\cdot, \mathbb{R}^d)$ is almost surely Lebesgue measure on $\mathbb{R}^d$. As $m \to \infty$, $\mathcal{T}_m(\mathbb{R}^d, \cdot)$ converges weak* to a Poisson point process of intensity 1 on $\mathbb{R}^d$. Finally, for any Borel set $A \subseteq \mathbb{R}^d$ with Lebesgue measure $\mathcal{L}(A) \in (0, \infty)$, we have

$$(11) \qquad \mathbf{E} \int\!\!\int \exp(c\|x-y\|^d) \mathcal{T}_m(dx, dy) \mathbf{1}[x \in A] \leq c' \mathcal{L}(A).$$

[To check (11) we first use invariance to deduce that the left-hand side must be a linear multiple of $\mathcal{L}(A)$, and then take $A$ to be a cube of volume $m$ to find the constant.]

By (11), the sequence $(\mathcal{T}_m)$ is tight, so let $\mathcal{T}$ be a weak* limit point, and note the following properties of $\mathcal{T}$. It is invariant, since invariant random transports form a weak* closed set. Clearly $\mathcal{T}(\cdot, \mathbb{R}^d)$ is Lebesgue measure on $\mathbb{R}^d$ almost surely. Writing $\Pi(\cdot) = \mathcal{T}(\mathbb{R}^d, \cdot)$, we see that $\Pi$ is a Poisson point process of intensity 1 on $\mathbb{R}^d$. And finally (11) holds with $\mathcal{T}$ in place of $\mathcal{T}_m$, since the set of random transports for which (11) holds is weak* closed.

Choose $\ell$ such that $1 - p = e^{-\ell^d}$, and for $z \in \mathbb{Z}^d$ define the unit cube $Q_z = [0, 1)^d + z$. Define a discrete configuration $\Gamma$ by

$$\Gamma(z) = 1 \wedge \Pi(\ell Q_z).$$

The choice of $\ell$ ensures that the law of $\Gamma$ is product measure with parameter $p$ on $\{0, 1\}^{\mathbb{Z}^d}$. Now define $\Theta$ by

$$\Theta_\Gamma(x, y) = \mathbf{E}(\mathcal{T}(\ell Q_x, \ell Q_y) | \Gamma).$$

It is elementary to check that $\Theta$ is a balanced transport rule for $\Gamma$, and (11) implies that it satisfies the required bound. $\square$

The following continuum analogue of Theorem 5 may be proved by applying Theorem 16 to the continuum transport given by

$$\Theta_\Pi(A, B) = \mathbf{E}(\mathcal{T}(A, B) | \Pi),$$

where $\mathcal{T}, \Pi$ are as in the above proof. Alternatively, it may be deduced from Theorem 5 by techniques similar to those used in [5], Section 4.

THEOREM 26. *Let $\Pi$ be a Poisson process of unit intensity on $\mathbb{R}^d$. If $d \geq 3$, then there exists a continuum extra head scheme for $\Pi$ satisfying*

$$\mathbf{E} \exp(C\|Y\|^d) < \infty$$

*for some $C = C(d) > 0$.*

**Open problems.**

(i) In the case of product measure on $\mathbb{Z}^d$ or a Poisson process on $\mathbb{R}^d$, what is the optimal tail behavior for *nonrandomized* extra schemes (or equivalently, for balanced allocation rules)?

(ii) What is the tail behavior of the extra head schemes (or allocation rules) constructed in Sections 4 and 5?

(iii) What is the optimal tail behavior of extra head schemes for product measure on other groups (e.g., for a free group with distance measured according to a Cayley graph)?

**Acknowledgments.** We thank Tom Liggett and Hermann Thorisson for introducing us to the problems, and for valuable discussions.

## REFERENCES

[1] BENJAMINI, I., LYONS, R., PERES, Y. and SCHRAMM, O. (1999). Group-invariant percolation on graphs. *Geom. Funct. Anal.* **9** 29–66. MR1675890
[2] GALE, D. and SHAPLEY, L. (1962). College admissions and stability of marriage. *Amer. Math. Monthly* **69** 9–15.
[3] HÄGGSTRÖM, O. (1997). Infinite clusters in dependent automorphism invariant percolation on trees. *Ann. Probab.* **25** 1423–1436. MR1457624
[4] HOFFMAN, C., HOLROYD, A. E. and PERES, Y. (2004). A stable marriage of Poisson and Lebesgue. Unpublished manuscript.
[5] HOLROYD, A. E. and LIGGETT, T. M. (2001). How to find an extra head: Optimal random shifts of Bernoulli and Poisson random fields. *Ann. Probab.* **29** 1405–1425. MR1880225
[6] HOLROYD, A. E. and PERES, Y. (2003). Trees and matchings from point processes. *Electron. Comm. Probab.* **8** 17–27. MR1961286
[7] KALLENBERG, O. (2002). *Foundations of Modern Probability. Probability and Its Applications*, 2nd ed. Springer, New York. MR1876169
[8] LIGGETT, T. M. (2002). Tagged particle distributions or how to choose a head at random. In *In and Out of Equilibrium* (V. Sidoravicious, ed.) 133–162. Birkhäuser, Boston. MR1901951
[9] MATTILA, P. (1995). *Geometry of Sets and Measures in Euclidean Spaces*. Cambridge Univ. Press. MR1333890
[10] MESHALKIN, L. D. (1959). A case of isomorphism of Bernoulli schemes. *Dokl. Akad. Nauk. SSSR* **128** 41–44. MR110782
[11] TALAGRAND, M. (1994). The transportation cost from the uniform measure to the empirical measure in dimension $\geq 3$. *Ann. Probab.* **22** 919–959. MR1288137
[12] THORISSON, H. (1995). On time- and cycle-stationarity. *Stochastic Process. Appl.* **55** 183–209. MR1313019
[13] THORISSON, H. (1996). Transforming random elements and shifting random fields. *Ann. Probab.* **24** 2057–2064. MR1415240

DEPARTMENT OF MATHEMATICS
UNIVERSITY OF BRITISH COLUMBIA
VANCOUVER, BRITISH COLUMBIA V6T 1Z2
CANADA
E-MAIL: holroyd@math.ubc.ca

DEPARTMENTS OF STATISTICS
AND MATHEMATICS
UNIVERSITY OF CALIFORNIA, BERKELEY
BERKELEY, CALIFORNIA 94720
USA
E-MAIL: peres@stat.berkeley.edu